\documentclass[12pt]{article}
\usepackage{amsfonts}
\usepackage{amsmath}
\usepackage{amssymb}

\input{tcilatex}
\input tcilatex
\begin{document}

\author{Refik Keskin$^{a}$ and Zafer \c{S}iar$^{b}$ \and $\ ^{a}$Sakarya
University, Mathematics Department, Sakarya/Turkey \and $^{b}$Bing\"{o}l
University, Mathematics Department, Bing\"{o}l/Turkey \and $^{a}$%
rkeskin@sakarya.edu.tr, $^{b}$zsiar@bingol.edu.tr}
\title{A Note on Terai's Conjecture Concerning the Exponential Diophantine
Equation $x^{2}+b^{y}=c^{z}$}
\maketitle

\begin{abstract}
Let $(a,b,c)$ be a primitive Pythagorean triple, i.e., $a^{2}+b^{2}=c^{2}$
with $\gcd (a,b,c)=1$, $a$ even and $b$ odd. Terai's conjecture says that
the Diophantine equation%
\begin{equation*}
x^{2}+b^{y}=c^{z}
\end{equation*}%
has only the positive integer solutions $(x,y,z)=(a,2,2)$. In this study, we
prove that Terai's conjecture is true when $b$ is a product of two primes
and $c\equiv 5(\func{mod}8)$.
\end{abstract}

\section{Introduction}

Let $(a,b,c)$ be a primitive Pythagorean triple, i.e., $a,b,c$ are positive
integers such that $a^{2}+b^{2}=c^{2}$ with $\gcd (a,b,c)=1$, $a$ even and $%
b $ odd. In \cite{Terai}, Terai conjecture that the Diophantine equation%
\begin{equation*}
x^{2}+b^{y}=c^{z}
\end{equation*}%
has only the positive integer solutions $(x,y,z)=(a,2,2)$ and proved that
the conjecture is true if $b$ and $c$ are primes with $b^{2}+1=2c.$ Terai's
conjecture is proved in the special cases. In \cite{Cao}, the authors proved
that the conjecture is true when $c\equiv 5(\func{mod}8)$ and $b$ or $c$ is
a prime power. In \cite{Yuan}, the authors showed that Terai's conjecture is
true if $b\equiv \pm 5(\func{mod}8)$ and $c$ is prime. Recently, in \cite%
{gok}, the authors have proved that Terai's conjecture is true if
\begin{equation*}
(a,b,c)=(2^{r+1}s,2^{2r}s^{2}-1,2^{2r}s^{2}+1),
\end{equation*}%
where $r,s\in
\mathbb{N}
,\gcd (2,s)=1,r\geq 2$, and$~s<2^{r-1}$. One can consult \cite{gok2} and
\cite{Le} for more information about Terai's conjecture.

It is well known that if $(a,b,c)$ is a primitive Pythagorean triple, then
there are positive integers $m$ and $n$ such that
\begin{equation*}
\gcd (m,n)=1,m>n,m\not\equiv n(\func{mod}%
2),a=2mn,b=m^{2}-n^{2},c=m^{2}+n^{2}.
\end{equation*}
In this study, we prove that Terai's conjecture is true when $b$ is a
product of two primes and $c\equiv 5(\func{mod}8)$. That is, we prove that

\begin{theorem}
\label{T1}Let $m,n\in
\mathbb{Z}
^{+}$ such that $m>n,~\gcd (m,n)=1$ and $m\not\equiv n(\func{mod}2).$ If $%
m^{2}+n^{2}\equiv 5(\func{mod}8),~m+n$ and $m-n$ are prime numbers, then the
only solution of the equation%
\begin{equation}
x^{2}+(m^{2}-n^{2})^{y}=(m^{2}+n^{2})^{z}  \label{3.1}
\end{equation}%
is $(x,y,z)$ $=(2mn,2,2).$
\end{theorem}

\section{Preliminaries}

In this section, we give some lemmas which will be useful in proving our
theorem.

\begin{lemma}
\label{L1}$\emph{\bigskip (\cite{Mor})}\mathbb{\ }$Let $x,y,z\in
\mathbb{Z}
,~z>0,~\gcd (x,y)=1$ and $2|y.$ Then every solution $(x,y,z)$ of the equation%
\begin{equation*}
x^{2}+y^{2}=z^{m}
\end{equation*}%
can be expressed as%
\begin{eqnarray*}
z &=&u^{2}+v^{2},~u,v\in
\mathbb{N}
,~\gcd (u,v)=1,~2|v, \\
x+iy &=&\lambda _{1}(u+\lambda _{2}iv)^{m},\lambda _{1},\lambda _{2}\in
\left\{ -1,1\right\} .
\end{eqnarray*}
\end{lemma}

In \cite{Zhu}, the authors proved the following result (see Corollary $1.2)$.

\begin{lemma}
\label{L2}Let $p$ and $n>3$ be primes. If $y$ is not the sum of two
consecutive squares, then the equation
\begin{equation*}
x^{2}+p^{2m}=2y^{n}
\end{equation*}%
has no solution $(x,y,m,n)$ in positive integers $x,y,m$ with $\gcd (x,y)=1$
$.$
\end{lemma}

The following lemma is proved by Cohn in \cite{Co} .

\begin{lemma}
\label{L3} The Diophantine equation $2z^{k}=y^{2}+1$ with $k>2$ has only the
solutions $y=z=1$ and $y=239,z=13,k=4.$
\end{lemma}

Thus, we can give the following result.

\begin{corollary}
\label{C}If $q\geq 3$ is odd, then the equation $y^{q}=a^{2}+(a+1)^{2}$ has
no positive integer solutions $y$ and $a.$
\end{corollary}

\section{The proof of Theorem \protect\ref{T1}}

Assume that $m^{2}+n^{2}\equiv 5(\func{mod}8),~\gcd (m,n)=1$, $m\not\equiv n(%
\func{mod}2)$, and $m+n=p,$ $m-n=q$ with prime integers $p,q.$ Let $(x,y,z)$
be a solution of the equation (\ref{3.1}). In this case,
\begin{equation*}
x^{2}\equiv -(m^{2}-n^{2})^{y}\equiv -(2n^{2})^{y}\left( \func{mod}%
m^{2}+n^{2}\right) .
\end{equation*}%
This implies that
\begin{equation*}
1=\left( \frac{-(2n^{2})^{y}}{m^{2}+n^{2}}\right) =\left( \frac{2}{%
m^{2}+n^{2}}\right) ^{y}=(-1)^{y}
\end{equation*}%
since $m^{2}+n^{2}\equiv 5(\func{mod}8),$ where the symbol $\left( \frac{%
\cdot }{\cdot }\right) $ denotes the Jacobi symbol. The last equality gives
that $y$ is even. Similarly, since
\begin{equation*}
x^{2}\equiv (m^{2}+n^{2})^{z}\equiv (2n^{2})^{z}\left( \func{mod}%
m^{2}-n^{2}\right)
\end{equation*}%
and $m^{2}-n^{2}=m^{2}+n^{2}-2n^{2}\equiv \pm 5(\func{mod}8),$ it follows
that
\begin{equation*}
1=\left( \frac{(2n^{2})^{z}}{m^{2}-n^{2}}\right) =\left( \frac{2}{m^{2}-n^{2}%
}\right) ^{z}=(-1)^{z}.
\end{equation*}%
Thus, $z$ is also even. So, let $y=2r$ and $z=2k$ for some positive integers
$r$ and $k.$ In this case, we have
\begin{equation}
x^{2}+(m^{2}-n^{2})^{2r}=(m^{2}+n^{2})^{2k}.  \label{3.2}
\end{equation}%
Since $(m,n)=1$ and $m\not\equiv n(\func{mod}2),$ it can be easily seen that
$\gcd (x,m^{2}-n^{2},m^{2}+n^{2})=1$ and so $(x,\left( m^{2}-n^{2}\right)
^{r},\left( m^{2}+n^{2}\right) ^{k})$ is a primitive Pythagorean triple.
Therefore, there exist positive integers $u$ and $v$ such that $%
x=2uv,~\left( m^{2}-n^{2}\right) ^{r}=u^{2}-v^{2},$ and $\left(
m^{2}+n^{2}\right) ^{k}=u^{2}+v^{2},$ where $u>v,~\gcd (u,v)=1$ and$%
~u\not\equiv v(\func{mod}2).$ The equalities $\left( m^{2}-n^{2}\right)
^{r}=u^{2}-v^{2}$ and $\left( m^{2}+n^{2}\right) ^{k}=u^{2}+v^{2}$ yield to
\begin{eqnarray*}
2u^{2} &=&\left( m^{2}-n^{2}\right) ^{r}+\left( m^{2}+n^{2}\right) ^{k} \\
2v^{2} &=&\left( m^{2}+n^{2}\right) ^{k}-\left( m^{2}-n^{2}\right) ^{r}.
\end{eqnarray*}%
From here, we get
\begin{equation*}
2u^{2}\equiv (m^{2}-n^{2})^{r}\equiv -(2n^{2})^{r}\left( \func{mod}%
m^{2}+n^{2}\right)
\end{equation*}%
and
\begin{equation*}
2v^{2}\equiv (m^{2}+n^{2})^{k}\equiv (2n^{2})^{k}\left( \func{mod}%
m^{2}-n^{2}\right) .
\end{equation*}%
These congruences, respectively, imply that
\begin{equation*}
-1=\left( \frac{2}{m^{2}+n^{2}}\right) =\left( \frac{-(2n^{2})^{r}}{%
m^{2}+n^{2}}\right) =-\left( \frac{2}{m^{2}+n^{2}}\right) ^{r}=(-1)^{r}
\end{equation*}%
and
\begin{equation*}
-1=\left( \frac{2}{m^{2}-n^{2}}\right) =\left( \frac{(2n^{2})^{k}}{%
m^{2}-n^{2}}\right) =\left( \frac{2}{m^{2}-n^{2}}\right) ^{k}=(-1)^{k}.
\end{equation*}%
Therefore, $r$ and $k$ are odd. On the other hand, we have
\begin{equation*}
u^{2}-v^{2}=(u-v)(u+v)=\left( m^{2}-n^{2}\right) ^{r}=p^{r}q^{r}.
\end{equation*}%
Since $\gcd (u-v,u+v)=1,$ it follows that either $u+v=p^{r},~u-v=q^{r}$ or $%
u-v=1$, $u+v=p^{r}q^{r}.$ Assume that $u-v=1$ and $u+v=p^{r}q^{r}.$ In this
case, we have
\begin{equation*}
2\left( m^{2}+n^{2}\right) ^{k}=2(u^{2}+v^{2})=(pq)^{2r}+1.
\end{equation*}%
Since $k$ is odd, this is impossible by Lemma \ref{L3} if $k>1.$ Therefore $%
k=1.$ Hence, we get
\begin{equation}
2\left( m^{2}+n^{2}\right) =(pq)^{2r}+1.  \label{3.12}
\end{equation}%
As $m+n=p,$ $m-n=q,$ we see that $2\left( m^{2}+n^{2}\right) =p^{2}+q^{2}.$
It follows that
\begin{equation}
(pq)^{2r}+1=p^{2}+q^{2},  \label{3.14}
\end{equation}%
which is impossible. Now assume that $u-v=q^{r}$ and$~u+v=p^{r}.$ Then, it
is seen that
\begin{equation}
p^{2r}+q^{2r}=2(u^{2}+v^{2})=2\left( m^{2}+n^{2}\right) ^{k}  \label{3.0}
\end{equation}%
and from here, we get%
\begin{equation*}
\left( \frac{p^{r}+q^{r}}{2}\right) ^{2}+\left( \frac{p^{r}-q^{r}}{2}\right)
^{2}=\left( m^{2}+n^{2}\right) ^{k}.
\end{equation*}%
Now, taking
\begin{equation*}
x=\frac{p^{r}+q^{r}}{2},~y=\frac{p^{r}-q^{r}}{2},z=m^{2}+n^{2}
\end{equation*}%
in Lemma \ref{L1} and considering the fact that $k$ is odd, we can write%
\begin{equation*}
m^{2}+n^{2}=z=c^{2}+d^{2},c,d\in \mathbb{Z},\gcd (c,d)=1,2|d,x+iy=(c+id)^{k}
\end{equation*}%
Thus, we can see that
\begin{equation*}
p^{r}=\frac{(c+di)^{k}+(c-di)^{k}}{2}+\frac{(c+di)^{k}-(c-di)^{k}}{2i}
\end{equation*}%
and
\begin{equation*}
q^{r}=\frac{(c+di)^{k}+(c-di)^{k}}{2}-\frac{(c+di)^{k}-(c-di)^{k}}{2i}.
\end{equation*}%
From the last equations, we get
\begin{equation}
2p^{r}=(c+di)^{k}(1-i)+(c-di)^{k}(1+i)  \label{3.4}
\end{equation}%
and
\begin{equation}
2q^{r}=(c+di)^{k}(1+i)+(c-di)^{k}(1-i).  \label{3.5}
\end{equation}%
Also, it can be seen that
\begin{eqnarray}
c+id &\equiv &-d+id\equiv d(-1+i)(\func{mod}c+d),  \label{3.6} \\
c-id &\equiv &-d-id\equiv -d(1+i)(\func{mod}c+d),
\end{eqnarray}%
\begin{equation}
c+id\equiv d+id\equiv d(1+i)(\func{mod}c-d),  \label{3.8}
\end{equation}%
and
\begin{equation}
c-id\equiv d-id\equiv d(1-i)(\func{mod}c-d).  \label{3.9}
\end{equation}%
Since $k$ is odd, it follows that $k=4t+1$ or $k=4t+3$ for some nonegative
integer $t.$ Let $k=4t+1.$ Then, using the congruences (\ref{3.6})-(\ref{3.9}%
), we obtain%
\begin{eqnarray*}
2p^{r} &=&(c+di)^{k}(1-i)+(c-di)^{k}(1+i) \\
&\equiv &d^{k}(-1+i)^{k}(1-i)+(-d)^{k}(1+i)^{k+1}(\func{mod}c+d) \\
&\equiv &(-d)^{k}\left( (1-i)^{k+1}+(1+i)^{k+1}\right) (\func{mod}c+d) \\
&\equiv &0(\func{mod}c+d)
\end{eqnarray*}%
and
\begin{eqnarray*}
2q^{r} &=&(c+di)^{k}(1+i)+(c-di)^{k}(1-i) \\
&\equiv &d^{k}(1+i)^{k}(1+i)+d^{k}(1-i)^{k}(1-i)(\func{mod}c-d) \\
&\equiv &d^{k}\left( (1+i)^{k+1}+(1-i)^{k+1}\right) (\func{mod}c-d) \\
&\equiv &0(\func{mod}c-d).
\end{eqnarray*}%
These imply that $c+d|2p^{r}$ and $c-d|2q^{r}.$ Since $c+d$ and $c-d$ are
odd, we see that $c+d|p^{r}$ and $c-d|q^{r}.~$Similarly, if $k=4t+3$, then
it can be seen that $c-d|p^{r}$ and $c+d|q^{r}$. Thus, we can write $%
\left\vert c+\epsilon d\right\vert =p^{t_{1}}$ and $\left\vert c-\epsilon
d\right\vert =q^{t_{2}}$ for some nonnegative integers $t_{1}$ and $t_{2}$,
where $\epsilon =\mp 1.$ Assume that $t_{2}=0$. If $t_{1}=0$, then we get $%
2=c^{2}+d^{2}$, which is impossible since
\begin{equation*}
c^{2}+d^{2}=m^{2}+n^{2}\geq 5.
\end{equation*}%
Let $t_{1}\geq 1$. Then
\begin{equation*}
p^{2}+q^{2}=2(m^{2}+n^{2})=2(c^{2}+d^{2})=\left\vert c-\epsilon d\right\vert
^{2}+\left\vert c+\epsilon d\right\vert ^{2}=1+\left\vert c+\epsilon
d\right\vert ^{2}=1+p^{2t_{1}},
\end{equation*}%
which is impossible if $t_{1}=1.$ If $t_{1}\geq 2$, then
\begin{equation*}
p^{2}+q^{2}=1+p^{2t_{1}}\geq 1+p^{4}>2p^{2},
\end{equation*}%
a contradiction. Let\ $t_{2}\geq 1$ and $t_{1}=0$. So we obtain%
\begin{equation*}
p^{2}+q^{2}=1+q^{2t_{2}}.
\end{equation*}%
It is seen that%
\begin{equation*}
(q^{t_{2}}-1)^{2}<q^{2t_{2}}-q^{2}+1=p^{2}<q^{2t_{2}},
\end{equation*}%
which implies that $q^{t_{2}}-1<p<q^{t_{2}}$. This is a contradiction. So, $%
t_{1}\geq 1$ and $t_{2}\geq 1$. Assume that $t_{1}>1$ or $t_{2}>1$. Then we
get
\begin{equation*}
p^{2}+q^{2}<p^{2t_{1}}+q^{2t_{2}}=\left\vert c+\epsilon d\right\vert
^{2}+\left\vert c+\epsilon d\right\vert
^{2}=2(c^{2}+d^{2})=2(m^{2}+n^{2})=p^{2}+q^{2},
\end{equation*}%
a contradiction. Therefore $t_{1}=t_{2}=1.$ This yields to $\left\vert
c+\epsilon d\right\vert =p$ and $\left\vert c-\epsilon d\right\vert =q$ and
so $m^{2}-n^{2}=pq=|c^{2}-d^{2}|$. As $m^{2}-n^{2}=|c^{2}-d^{2}|$ and $%
m^{2}+n^{2}=c^{2}+d^{2}$, it can be seen that\ $cd=\pm mn$. Now, let $\alpha
=(c+id)$ and $\beta =(c-id).$ Multiplying the equations (\ref{3.4}) and (\ref%
{3.5}) side to side, we obtain%
\begin{equation}
4p^{r}q^{r}=4(m^{2}-n^{2})^{r}=2(\alpha ^{2k}+\beta ^{2k}).  \label{3.11}
\end{equation}%
Now assume that $k>1$ and $r>1.$ Then, since $m^{2}-n^{2}=|c^{2}-d^{2}|,$ it
follows that
\begin{eqnarray*}
2(m^{2}-n^{2})^{r} &=&\alpha ^{2k}+\beta ^{2k}=(c+id)^{2k}+(c-id)^{2k} \\
&=&\left( (c+id)^{2}+(c-id)^{2}\right) \cdot A \\
&=&2\left( c^{2}-d^{2}\right) \cdot A=\mp (m^{2}-n^{2})A,
\end{eqnarray*}%
which implies that $A=\pm (m^{2}-n^{2})^{r-1}$, where%
\begin{equation*}
A=\dsum\limits_{j=0}^{k-1}\left( \alpha ^{2}\right) ^{k-1-j}\left( -\beta
^{2}\right) ^{j}.
\end{equation*}%
Note that $\alpha ^{2}+\beta ^{2}=2(c^{2}-d^{2}).$ On the other hand, since $%
2\left( m^{2}-n^{2}\right) =2\left\vert c^{2}-d^{2}\right\vert =\pm \left(
\alpha ^{2}+\beta ^{2}\right) ,$ it is seen that
\begin{equation*}
\alpha ^{2}\equiv -\beta ^{2}(\func{mod}m^{2}-n^{2}).
\end{equation*}%
Therefore,
\begin{eqnarray*}
A &=&\dsum\limits_{j=0}^{k-1}\left( \alpha ^{2}\right) ^{k-1-j}\left( -\beta
^{2}\right) ^{j}\equiv \dsum\limits_{j=0}^{k-1}\left( -\beta ^{2}\right)
^{k-1-j}\left( -\beta ^{2}\right) ^{j}(\func{mod}m^{2}-n^{2}) \\
&\equiv &\dsum\limits_{j=0}^{k-1}\left( -\beta ^{2}\right) ^{k-1}(\func{mod}%
m^{2}-n^{2}) \\
&\equiv &k\cdot \left( \beta ^{2}\right) ^{k-1}(\func{mod}m^{2}-n^{2}).
\end{eqnarray*}%
Also, since $\beta ^{2}=c^{2}-d^{2}-2cdi$ and $cd=\pm mn$, we see that
\begin{equation*}
\beta ^{2}\equiv \mp 2mni(\func{mod}m^{2}-n^{2})
\end{equation*}%
and thus, $\left( \beta ^{2}\right) ^{k-1}\equiv 2^{k-1}\cdot (mn)^{k-1}(%
\func{mod}m^{2}-n^{2}),$ where we used the fact that that $k-1$ is even.
Hence,
\begin{equation*}
A\equiv k\cdot 2^{k-1}\cdot (mn)^{k-1}(\func{mod}m^{2}-n^{2}).
\end{equation*}%
As $r>1$, and $A=\pm (m^{2}-n^{2})^{r-1}$, it is seen that $m^{2}-n^{2}|A.$
So it follows that $m^{2}-n^{2}|k\cdot (mn)^{k-1}$ since $m^{2}-n^{2}$ is
odd. It can be easily seen that
\begin{equation*}
\gcd ((mn)^{k-1},m^{2}-n^{2})=1.
\end{equation*}%
Thus, we get $m^{2}-n^{2}|k.$ Then,
\begin{equation*}
k=(m^{2}-n^{2})a=pqa
\end{equation*}%
for some odd positive integer $a.$ Let rearrange the equation (\ref{3.0}) as
\begin{equation}
2\left( \left( m^{2}+n^{2}\right) ^{qa}\right) ^{p}=p^{2r}+q^{2r}.
\label{3.10}
\end{equation}%
Since $q\geq 3$, it is seen that $\left( m^{2}+n^{2}\right) ^{qa}$ is not
the sum of two consecutive squares by Corollary \ref{C}. Then by Lemma \ref%
{L2}, the equation (\ref{3.10}) has no solutions since $p>3.$

Consequently, $k=1$ or $r=1.$ If $k=1$, then (\ref{3.0}) implies that
\begin{equation}
p^{2r}+q^{2r}=2(u^{2}+v^{2})=2\left( m^{2}+n^{2}\right)
^{k}=2(m^{2}+n^{2})=p^{2}+q^{2}.  \label{3.13}
\end{equation}%
This yields to $r=1.$ If $r=1$, then (\ref{3.0}) shows that
\begin{equation*}
2\left( m^{2}+n^{2}\right)
=p^{2}+q^{2}=p^{2r}+q^{2r}=2(u^{2}+v^{2})=2(m^{2}+n^{2})^{k}
\end{equation*}%
and therefore $k=1.$ Thus $y=2r=2,z=2k=2$ and so the only solution is $%
(x,y,z)$ $=(2mn,2,2)$.%
\endproof%

\end{document}